\documentclass[conference]{IEEEtran}
\IEEEoverridecommandlockouts
\usepackage{cite}
\usepackage{amsmath,amssymb,amsfonts}
\usepackage{algorithmic}
\usepackage{graphicx}
\usepackage{textcomp}
\usepackage{xcolor}
\usepackage{hyperref}
\usepackage{graphicx}
\usepackage[caption=false]{subfig}
\allowdisplaybreaks
\makeatletter
\newcommand{\linebreakand}{%
  \end{@IEEEauthorhalign}
  \hfill\mbox{}\par
  \mbox{}\hfill\begin{@IEEEauthorhalign}
}
\makeatother

\def\BibTeX{{\rm B\kern-.05em{\sc i\kern-.025em b}\kern-.08em
    T\kern-.1667em\lower.7ex\hbox{E}\kern-.125emX}}

\begin{document}

\title{Parallelizing Explicit and Implicit Extrapolation Methods for Ordinary Differential Equations
}

\author{
\IEEEauthorblockN{Utkarsh}
\IEEEauthorblockA{Massachusetts Institute of Technology\\Julia Computing\\
Cambridge, Massachusetts, USA \\
\href{mailto:utkarsh5@mit.edu}{utkarsh5@mit.edu}}
\and
\IEEEauthorblockN{Chris Elrod}
\IEEEauthorblockA{Julia Computing
\\Cambridge, Massachusetts, USA\\
\href{mailto:chris.elrod@juliacomputing.com}{chris.elrod@juliacomputing.com}}
\and
\IEEEauthorblockN{Yingbo Ma}
\IEEEauthorblockA{Julia Computing
\\Cambridge, Massachusetts, USA\\
\href{mailto:yingbo.ma@juliacomputing.com}{yingbo.ma@juliacomputing.com}}
\linebreakand
\IEEEauthorblockN{Konstantin Althaus}
\IEEEauthorblockA{Technical University of Munich
\\Munich, Germany\\
\href{mailto:konstantin.althaus@tum.de}{konstantin.althaus@tum.de}}
\and
\IEEEauthorblockN{Christopher Rackauckas}
\IEEEauthorblockA{Julia Computing\\Massachusetts Institute of Technology\\Pumas AI Inc.
\\Cambridge, Massachusetts, USA\\
\href{mailto:chris.rackauckas@juliacomputing.com}{chris.rackauckas@juliacomputing.com}}
}

\maketitle

\begin{abstract}
Numerically solving ordinary differential equations (ODEs) is a naturally serial process and as a result the vast majority of ODE solver software are serial. In this manuscript we developed a set of parallelized ODE solvers using extrapolation methods which exploit ``parallelism within the method'' so that arbitrary user ODEs can be parallelized. We describe the specific choices made in the implementation of the explicit and implicit extrapolation methods which allow for generating low overhead static schedules to then exploit with optimized multi-threaded implementations. We demonstrate that while the multi-threading gives a noticeable acceleration on both explicit and implicit problems, the explicit parallel extrapolation methods gave no significant improvement over state-of-the-art even with a multi-threading advantage against current optimized high order Runge-Kutta tableaus. However, we demonstrate that the implicit parallel extrapolation methods are able to achieve state-of-the-art performance (2x-4x) on standard multicore x86 CPUs for systems of $<200$ stiff ODEs solved at low tolerance, a typical setup for a vast majority of users of high level language equation solver suites. The resulting method is distributed as the first widely available open source software for within-method parallel acceleration targeting typical modest compute architectures.
\end{abstract}


\section{Introduction}



The numerical approximation of ordinary differential equations (ODEs) is a naturally serial time stepping process, meaning that methods for parallelizing the solution of such ODEs requires either tricks or alternative solvers in order to exploit parallelism. The most common domain, and one of the most common applications in scientific computing, is the parallel solution of large-scale systems of ODEs (systems of millions or more equations) which arise from the semi-discretization of partial differential equations. For such cases, the size of the system allows for parallel efficiency to be achieved by parallelizing some of the most expensive computations, such as implicit parallelism of BLAS \cite{dongarra1990set}/LAPACK \cite{anderson1999lapack}, sparse linear solvers, or preconditioner computations, on compute clusters and with heterogeneous GPU+CPU compute. Commonly used open source softwares such as Sundials \cite{Hindmarsh_2005} and DifferentialEquations.jl \cite{rackauckas2017differentialequations} help users automate the parallelism in such cases. However, for this manuscript we look in the opposite direction at parallelism for small systems of ODEs ($<200$).

Parallel computing on small systems of equations comes with a separate set of challenges, namely that computations can easily become dominated by overhead. For this reason, most software for accelerating the parallel solution of ODEs focus on the parallelization of generating ensembles of solutions, i.e. generating the solution to a small system of ODEs over a set of parameters or initial conditions. Examples of this include the ensembles interface in DifferentialEquations.jl \cite{rackauckas2017differentialequations}. However, in many cases the ensemble form may not be easy to exploit. One key example is in parameter estimation or model calibration routines which are typically done with some gradient-based optimization technique and requires one realization of the model at a given parameter before advancing. Given questions on the Julia \cite{bezanson2017julia} Discourse demonstrate that the vast majority of users both typically solve $<200$ ODE systems and have multi-core compute available (e.g. Core i5/i7 chips with 2-16 processors), investigated whether the most standard ODE solve cases could benefit from some form of parallelism.

Prior research has developed a large amount of theory around potential solvers for such problems which fall into the categories of parallel-in-time (PIT) methods and within-step parallelism methods \cite{burrage1995parallel}. Parallel-in-time methods address the issue by effectively stepping the ODE at multiple time points simultaneously, similar to the ensemble approach, and impose a nonlinear system constraint to relax the initial conditions of the future points to arrive at a sufficiently smooth solution. An open source software for PIT methods, XBraid, does exist but focuses on large compute hardware \cite{xbraid-package}. Documents from the XBraid tutorials \cite{xbraid-tutorial} estimate current PIT methods outperform classical solvers at $\sim 256$ cores\footnote{Private correspondences and tests with DifferentialEquations.jl also confirm a similar cutoff point with PIT methods}, making them impractical for standard compute hardware. Thus while PIT methods are a promising direction for accelerating ODE solving for the exascale computers targeted by the XBraid project \cite{xbraid-package}, these methods do not suffice for the $\leq 16$ core nature of standard consumer computing devices.

However, promising prior results exist within the within-step parallelism research. Nowak 1998 parallelized the implicit extrapolation methods and showed that the parallelization did improve performance \cite{nowak1998parallel}, though its tests did not compare against multithreaded (sparse) factorizations like is seen in modern suites such as UMFPACK \cite{davis2004algorithm_umfpack} in SuiteSparse \cite{davis2014suitesparse}, no comparison was made against optimized solver softwares such as CVODE \cite{Hindmarsh_2005} or LSODA \cite{petzold1983automatic}, and no open source implementation exists for this method. Additionally, Ketcheson (2014) \cite{ketcheson2014comparison} demonstrated that within-step parallel extrapolation methods for non-stiff ODEs can outperform dop853, an efficient high-order Runge-Kutta method \cite{hairer1993solving} . Two notable issues with this study are that (a) no optimized open source software exists for this method for further investigation, (b) the study did not investigate whether such techniques could be useful for stiff ODEs, which tend to be the most common type of ODEs for many applications from biology and chemistry to engineering models. These results highlight that it may be possible to achieve state-of-the-art (SOA) performance by exploiting parallelism, though no software is readily available.
 
In this manuscript we build off of the prior work in within-method parallel extrapolation to build an open source software in Julia \cite{bezanson2017julia} targeting $<200$ ODE systems on standard compute architectures. We demonstrate that this new solver is the most efficient solver for this class of equations when high accuracy is necessary, and demonstrate this with benchmarks that include many methods from DifferentialEquations.jl \cite{rackauckas2017differentialequations}, SUNDIALS \cite{Hindmarsh_2005} (CVODE), LSODA \cite{petzold1983automatic}, and more. To the best of our knowledge, this is the first demonstration that a within-step parallel solver can be the most efficient method for this typical ODE and hardware combination.

\section{Extrapolation Methods}

Our exposition of extrapolation methods largely follows that of Hairer I \cite{hairer1993solving} and II \cite{hairer1991solving}, and \cite{althaus_2018}. Take the ODE:
\begin{equation}
    u^{\prime} = f(u,p,t),
\end{equation}
with some known initial condition $u(t_0) = u_0$ over a time span $t \in (t_0, t_f)$. The extrapolation methods are variable-order and variable-step methods which generate higher precision approximations of ODE solutions computed at different time step-sizes. For the $N^{th}$ current order of the algorithm, we generate $f(u,\ p,\ t+dt)$ at the current time-step $dt$ for each order from 1 to $N$. For the approximation at $t+dt$, the algorithm chooses a subdividing sequence which discretizes into further fixed smaller steps between $t \ \text{and} \ t+dt$. Choosing a sequence of the form:
\begin{equation}\label{eq:1}
    n_1 < n_2 < n_3 < n_4 < \cdots < n_{N}
\end{equation}
Generates internal step-sizes $h_1 > h_2 > h_3 > h_4 > \cdots > h_{N}$ by $h_i = \frac{dt}{n_i}$. The subdividing sequences vary with order, having smaller time-steps in higher orders for finer resolution. The algorithm chosen for this is suitably an efficient implicit/explicit method between $t \ \text{and} \ t+dt$ is of $p^{th}$ order. 

The tabulation of $1^{st}$ $N^{th}$ calculations generate starting-stage of extrapolation computation, denoted by:
\begin{align}\label{eq:3}
    \begin{split}
        {}&u_{h_i}(t+dt) = T_{i,1}
        \\
        {}&i = j,j-1,j-2,\cdots,j-k+1.
        \\
    \end{split}
\end{align}
%
%
%
Extrapolation methods use the interpolating polynomial
\begin{equation}\label{eq:5}
    p(h) = \Tilde{u} - e_ph^p - e_{p+1}h^{p+1} - \cdots - e_{p+k-2}h^{p+k-2}
\end{equation}
such that $p(h_i) = T_{i, 1}$ to obtain a higher order approximation by extrapolating the step size $h$ to $0$.
Concretely, we define $T_{j,k} := p(0) = \Tilde{u} = u(t + dt) + \mathcal{O}(dt^{p + k})$. In order to find $\Tilde{u}$, we can solve the linear system of $k$ variables and $k$ equations formed by equations \eqref{eq:3} and \eqref{eq:5}. This conveniently generates an array of approximations with different orders which allows simple estimates of local error and order variability techniques:
\[    
\begin{matrix}
    T_{1,1} &  & \\
    T_{2,1} & T_{2,2} & \\
    T_{3,1} & T_{3,2} & T_{3,3} \\
    \cdots  & \cdots & \cdots & \cdots \\
\end{matrix}
\] 
Aitken-Neville's algorithm uses Lagrange polynomial interpolation formulae \cite{aitken1932interpolation,neville1934iterative} to make: \eqref{eq:5}:
\begin{equation}\label{eq:8}
    T_{j,k+1} = T_{j,k} + \frac{T_{j,k}-T_{j-1,k}}{\frac{n_j}{n_{j-k}}-1}
\end{equation}
Finally, the $N^{th}$ order is $u(t+dt) = T_{N,N}$.

As shown in \eqref{eq:1}, the subdividing sequence should be positive and strictly increasing. Common choices are:
\begin{enumerate}
    \item Harmonic \cite{deuflhard1983order}: $n = 1, 2, 3, 4, 5, 6, 7, 8 \hdots$
    \item Romberg \cite{romberg1955vereinfachte}: $n = 1, 2, 4, 8, 16, 32, 64, 128, 256 \hdots$
    \item Bulrisch \cite{romberg1955vereinfachte}, \cite{bulirsch1966numerical}: $n = 1, 2, 3, 4, 6, 8, 12, 16 \hdots$
\end{enumerate}
The "Harmonic" sequence generates the most efficient load balancing and utilization of multi-threading in parallel computing which is discussed more in Section III.

\subsection{Explicit Methods}

\subsubsection{Extrapolation Midpoint Deuflhard and HairerWanner}

Both the algorithms use explicit midpoint method for internal step-size computations. The representation used here is in the two-step form, which makes the algorithm symmetric and has even powers in asymptotic expansion \cite{stetter1970symmetric}
\begin{align}
    \begin{split}
        u_{h_i}(t_0) &= u_0
        \\
        u_{h_i}(t_1) &= u_0 + h_i f(u_0,p,t_0)
        \\
        u_{h_i}(t_n) &= u_{h_i}(t_{n-2}) + 2h_if(u_{h_i}(t_{n-1}),p,t_{n-1})
    \end{split}
\end{align}
The difference between them arises from the step-sizing controllers. The ExtrapolationMidpointDeufhard is based on the Deuflhard's DIFEX1 \cite{deuflhard1983order} adaptivity behaviour and ExtrapolationMidpointHairerWanner is based on Hairer's ODEX \cite{hairer1993solving} adaptivity behaviour.

The extrapolation is performed using barycentric formula which based on the lagrange barycentric interpolation \cite{berrut2004barycentric}. The interpolation polynomial is given by:
\begin{align}
    \begin{split}
        w_j &= \prod_{i=1:N+1,i\neq j}\frac{1}{n_j^{-2} - n_i^{-2}}
        \\
        \rho(h) &= \prod_{i=1:N+1}{h - n_i^{-2}}
        \\
        p(h) &= \rho(h)\sum^{N+1}_{j = 1}\frac{w_j}{h-{n_j}^{-2}}T_{j,1}
    \end{split}
\end{align}
Extrapolating the limit $h\rightarrow0$, we get:
\begin{align}\label{eq:11}
    \begin{split}
        u(t+dt) &= \rho(0)\sum^{N+1}_{j = 1}\frac{w_j}{-{n_j}^{-2}}T_{j,1}
    \end{split}
\end{align}
Where $w_j$ are the extrapolation weights, $\rho(0)$ are the extrapolation coefficients and $n_j$ denotes the subdividing sequence.

The choice for baryentric formula instead of Aitken-Neville is mainly due to reduced computation cost of ODE solution at each time-step \cite{althaus_2018}. The extrapolation weights $w_j$ and coefficients $\rho(0)$ can be easily computed and stored as tableau's. The yields the computation cost to be $O(N_{max}^2)$ for the precomputation where $N_{max}$ is the maximal order than can be achieved by the method.The extrapolation method of order $N$ generates a method of error in order of $2(N+1)$ \cite{hairer1993solving}, \cite{althaus_2018}.

It can be analysed that computational cost of with Aitken-Neville for $T_{N,N}$ is $\mathcal{O}(N^2d)$ ($d(1+2+\cdots+N) = \mathcal{O}(N^2d)$), where $N$ is the extrapolation order and d is the dimension of the $u$ \cite{althaus_2018}. 
The computation cost of \eqref{eq:11} is simply \textbf{$\mathcal{O}(Nd)$} (a linear combination all $T_{j,1}$ across d dimensions) \cite{althaus_2018}.
\bigbreak
\textbf{Numerical Stability and Analysis}: In comparison The numerical performance follows similar behaviour as that of Aitken-Neville in the subdividing sequences of Romberg and Bulrisch, where the absolute error decreases as extrapolation order $q$ increases \cite{webb2012stability}. However, in the harmonic sequence, the numerical stability in both cases of extrapolation remains up-to as much as up-to 15 order and then it diverges \cite{althaus_2018}.

\subsection{Implicit Methods}

We will consider implicit methods as basis for internal step-size calculation using subdividing sequences in extrapolation methods. They are widely used for solving stiff ODEs (maybe write more why implicit methods are suited for it).

\subsubsection{Implicit Euler Extrapolation}

The ImplicitEulerExtrapolation uses the Linearly-Implicit Euler Method for internal step-sizing. Mathematically:
\begin{align}
    \begin{split}
        (I - h_iJ)(u_{h_i}(t_{n+1}) - u_{h_i}(t_n)) = hf(u_{h_i}(t_n),p,t_n)
        \\
        \therefore u_{h_i}(t_{n+1}) = u_{h_i}(t_{n}) + (I - h_iJ)^{-1}hf(u_{h_i}(t_{n+1}),p,t_n)
    \end{split}
\end{align}
where $I$ and $J \approx \frac{\partial f(u_{h_i},p,t)}{\partial u_{h_i}} $ is the $\mathbb{R}^{d\times d}$ identity and jacobian matrix respectively. Clearly, the method is non-symmetrical and consequently would have non-even powers of $h$ global error expansion. The methods are $A(\alpha)$ stable with $\alpha \approx 90^o$ \cite{hairer1991solving}. The extrapolation scheme used is Aitken-Neville \eqref{eq:8}.

\subsubsection{Implicit Euler Barycentric Extrapolation}

We experimented with barycentric formulas \cite{berrut2004barycentric} to replace extrapolation algorithm in ImplicitEulerExtrapolaton. Since there are no even powers in global error expansion, the barycentric formula \cite{berrut2004barycentric} needs to changed as follows:
\begin{align}
    \begin{split}
        w_j &= \prod_{i=1:N+1,i\neq j}\frac{1}{n_j^{-1} - n_i^{-1}}
        \\
        \rho(h) &= \prod_{i=1:N+1}{h - n_i^{-1}}
        \\
        p(h) &= \rho(h)\sum^{N+1}_{j = 1}\frac{w_j}{h-{n_j}^{-1}}T_{j,1}
    \end{split}
\end{align}
Extrapolating the limit $h\rightarrow0$, we get:
\begin{align}
    \begin{split}
        u(t+dt) &= \rho(0)\sum^{N+1}_{j = 1}\frac{w_j}{-{n_j}^{-1}}T_{j,1}
    \end{split}
\end{align}
Consequently, the extrapolation method of order $N$ generates a method of error in order of $N+1$.


\subsubsection{Implicit Hairer Wanner Extrapolation}

For implicit extrapolation, symmetric methods would be beneficial to provide higher order approximations. The naive suitable candidate is the trapezoidal rule, but the resulting method is not stiffly stable and hence undesirable for solving stiff ODEs \cite{dahlquist1963special}, \cite{hairer1993solving}.

G. Bader and P. Deuflhard \cite{bader1983semi}, \cite{deuflhard1983order} developed the linearly-implicit midpoint rule with Gragg's smoothing \cite{gragg1964repeated}, \cite{Lindberg_1971} as extension of popular Gragg-Bulirsch-Stoer (GBS) \cite{bulirsch1966numerical} method. The algorithm is given by:
\begin{align}
    \begin{split}
        (I - h_iJ)(u_{h_i}(t_{1}) - u_0) &= h_if(u_0,p,t_0)
        \\
        (u_{h_i}(t_{n+1}) - u_{h_i}(t_n)) &= (u_{h_i}(t_n) - u_{h_i}(t_{n-1}))\\ 
        & + 2(I - h_iJ)^{-1}h_if(u_{h_i}(t_n),p,t_n)\\
        & - 2(I - h_iJ)^{-1}(u_{h_i}(t_n) - u_{h_i}(t_{n-1}))
    \end{split}
\end{align}
Followed by Gragg's Smoothing \cite{gragg1964repeated}, \cite{Lindberg_1971}:
\begin{align}
    \begin{split}
        t &= t_0 + 2n{h_i}\\
        T_{j,1} = S_{h_i}(t) &= \frac{u_{h_i}({2n+1}) + u_{h_i}({2n-1})}{2}\\
    \end{split}
\end{align}
The Gragg's smoothing \cite{gragg1964repeated}, \cite{Lindberg_1971} implementation in OrdinaryDiffEq.jl \footnote{\url{https://github.com/SciML/OrdinaryDiffEq.jl/pull/1212}} improved the stability of algorithm resulting in lower $f$ evaluations, time steps and linear solves. Consequently, it increased the accuracy and speed of the solvers.

The extrapolation scheme used is similar to extrapolation midpoint methods, namely the barycentric formula \cite{berrut2004barycentric} \eqref{eq:11}. One of the trade-offs for the increased stability from Gragg's smoothing \cite{gragg1964repeated}, \cite{Lindberg_1971} is the reduced order of error for the extrapolation method. As stated earlier, the barycentic formula provides and error of order $2(N+1)$ which gets reduced to $2(N+1)-1 = 2N + 1$ in the case of Implicit Hairer Wanner Extrapolation \cite{hairer1993solving}, \cite{bader1983semi}. The algorithm requires to use even subdividing sequence and hence we are using multiples of 4 of the common sequences in the implementation.
Furthermore, the convergence tests for these methods\footnote{\url{https://github.com/SciML/OrdinaryDiffEq.jl/blob/master/test/algconvergence/ode_extrapolation_tests.jl}} are written in the test suite of OrdinaryDiffEq.jl.

\subsection{Adaptive time-stepping and order of algorithms}

The methods follow a comprehensive step-sizing and order selection strategy. The adaptive time-stepping is mostly followed on the lines of Hairer's ODEX \cite{hairer1993solving} adaptivity behaviour. The error at $k^{th}$ order is expressed as:
\begin{align}
    \begin{split}
        err_k = \lVert T_{k,k-1} - T_{k,k} \rVert
    \end{split}
\end{align}
and the scaled error is represented as:
\begin{align}
    \begin{split}
        err_{scaled} = \frac{err_k}{abstol + max(u,u_{prev})reltol}
    \end{split}
\end{align}
Finally, the $h_{optimal}$ is calculated using the standard-controller:
\begin{align}
    \begin{split}
        &q = max(q_{min},min(\frac{err_{scaled}^{alg\_order+1}}{\gamma},q_{max}))
        \\
        &h_{optimal} = \frac{h}{q}
    \end{split}
\end{align}
Where $\gamma$ is the safety factor, $q_{min}$ is the minimum scaling and $q_{max}$ is the maximum scaling. For the optimal order selection, the computation relies on the "work calculation" which is stage-number $A_k$ per step-size $h_k$,  $\frac{A_k}{h_k}$. 

The stage-number includes the no. of $f(u,p,t)$ RHS function evaluations, jacobian matrix evaluation, and forward and backward substitutions. These are pre-computed according the subdividing sequence and stored as a cache. More details about these calculations can be found here.\footnote{\url{https://github.com/utkarsh530/DiffEqBenchmarks/blob/master/Extrapolation_Methods/Extrapolation_Methods.pdf}}
These work calculations form the basis of order-selection. The order-selection is restricted between the window $(k-1,k+1)$ and appropriate conditions are passed to convergence monitor. Briefly describing, it checks the convergence in the window and subsequently accepts the order which has the most significant work reduction during increasing/decreasing the order and $error < 1$. Detailed description can be found at order and step-size control in \cite{hairer1993solving}.

\section{Parallelization of the Algorithm}

\subsection{Choosing Subdividing Sequences for Static Load Balancing}

The exploit of parallelism arises from the computation of $T_{j,1}$. $k$ evaluations of $T_{j,1}$ are done to find the solution $u(t+dt)$. These evaluations are independent and are thus parallelizable. However, because of the small systems being investigated, dynamic load balancing overheads are much too high to be relied upon and thus we tailored the parallelism implementation to make use of a static load balancing scheme as follows. Each $T_{j,1}$ requires $n_j$ calls to the function $f$, and if the method is implicit then there is an additional $n_j$ back substitutions and one LU factorization. While for sufficiently large systems the LU-factorization is the dominant cost due to its $O(n^3)$ growth for an $n$ ODE system, in the $<200$ ODE regime with our BLAS \cite{dongarra1990set}/LAPACK \cite{anderson1999lapack} implementations\footnote{RecursiveFactorization.jl notably outperforms OpenBLAS in this regime.} we find that this cost is close enough to the cost of a back substitution that we can roughly assume $T_{j,1}$ requires $n_j$ ``work units'' in both the implicit and explicit cases method cases. 

If we consider the multiples of harmonic subdividing sequence is 2, 4, 6, 8, 10, 12, ... then the computation of $T_{j,1}$ needs $2j$ work. Because $T_{k,k}$ is a numerical approximation of order $p = 2k$ we can simply give each $T_{j,1}$ a processor if $k$ processors are available. The tasks would not finish simultaneously due to different amount of work. With Multiple Instruction and Multiple Data Stream (MIMD) processors, multiple computations can be loaded to one processor to calculate say $T_{1,1}$ and $T_{k-1,1}$ on a single processor, resulting in $2 + 2(k-1) = 2k$ function evaluations . Consequently, each processor would contain $T_{j,1}$ and $T_{k-j,1}$ computations and we require $\frac{k}{2}$ processors. This is an effective load-balancing and all processes would finish simultaneously due to an approximately constant amount of work per chunk.

Julia provides compose-able task-based parallelism. To take advantage of the regularity of our work and this optimized static schedule, we avoid some of the overheads associated with this model such as task-creation and scheduling by using Polyester.jl's \cite{ElrodPolyester} \texttt{@batch} macro, which allows us to run our computation on long-lived tasks that can remain active with a spin-lock between workloads, thereby avoiding the need for rescheduling. We note that this load-balancing is not as simple in sequences of Romberg and Bulirsch, and thus in those cases we multi-thread computations over constant $p$ processors only. For this reason, the homonic subdividing sequence is preferred for performance in our implementation.

\subsection{Parallelization of the LU Factorization}

One of the most important ways multi-threading is commonly used within implicit solvers is within the LU factorization. The LU factorization is multi-threaded due to expensive $O(N^3)$ complexity. For Jacobian matrices smaller than $100 \times 100$, commonly used BLAS/LAPACK implementations such as OpenBLAS, MKL, and RecursiveFactorization.jl are unable to achieve good parallel performance. However, the implicit extrapolation methods require solving linear systems with respect to $I - h_i J$ for each sub-time step $h_i$. Thus for the within-method parallel implicit extrapolation implementation we manually disabled the internal LU factorization multi-threading and parallelized the LU factorization step by multi-threading the computation of this ensemble of LU factorizations. If the chosen order is sufficiently then a large number of $h_i$ would be chosen and thus parallel efficiency would be achieved even for small matrices. In the results this will be seen as the lower tolerance solutions naturally require higher order methods, and this is the regime where the implicit extrapolation methods demonstrably become SOA.

\section{Benchmark Results}

The benchmarks we use work-precision diagrams to allow simultaneous comparisons between speed and accuracy. Accuracy is computed against reference solutions at $10^{-14}$ tolerance. All benchmarks were computed on a AMD EPYC 7513 32-Core processor ran at 8 threads\footnote{While more threads were available, we did not find more threads to be beneficial to accelerating these computations.}.

\subsection{Establishing Implementation Efficiency}

\begin{figure}[!t]
\centering
\includegraphics[width=\linewidth]{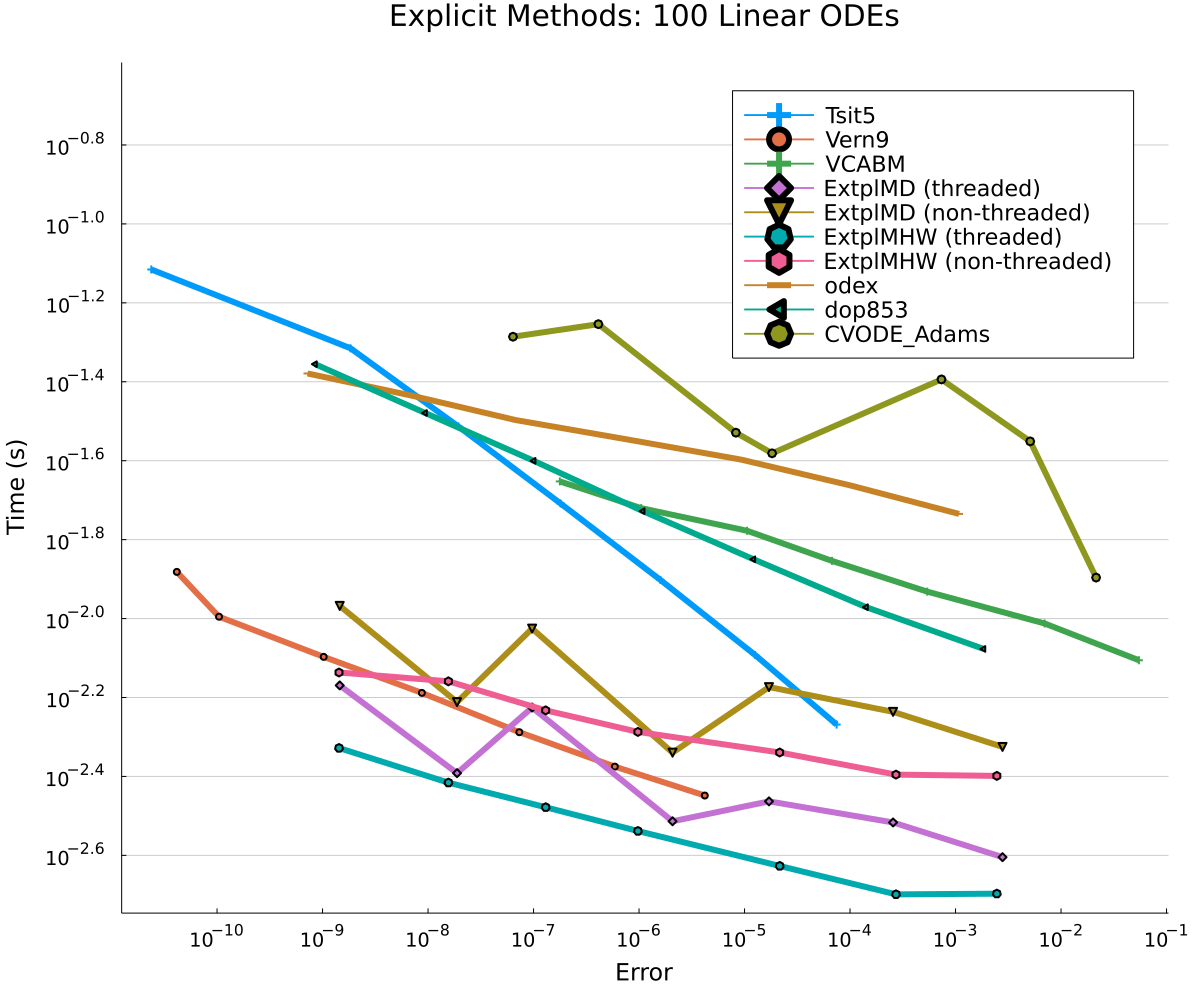}
\vspace{-25px}
\caption{Benchmark on the 100 linear ODE problem [\ref{sssec:linear}].}
\label{fig:linearcomparison}
\vspace{-10px}
\end{figure}

\begin{figure}[!t]
\centering
{\includegraphics[width=\linewidth]{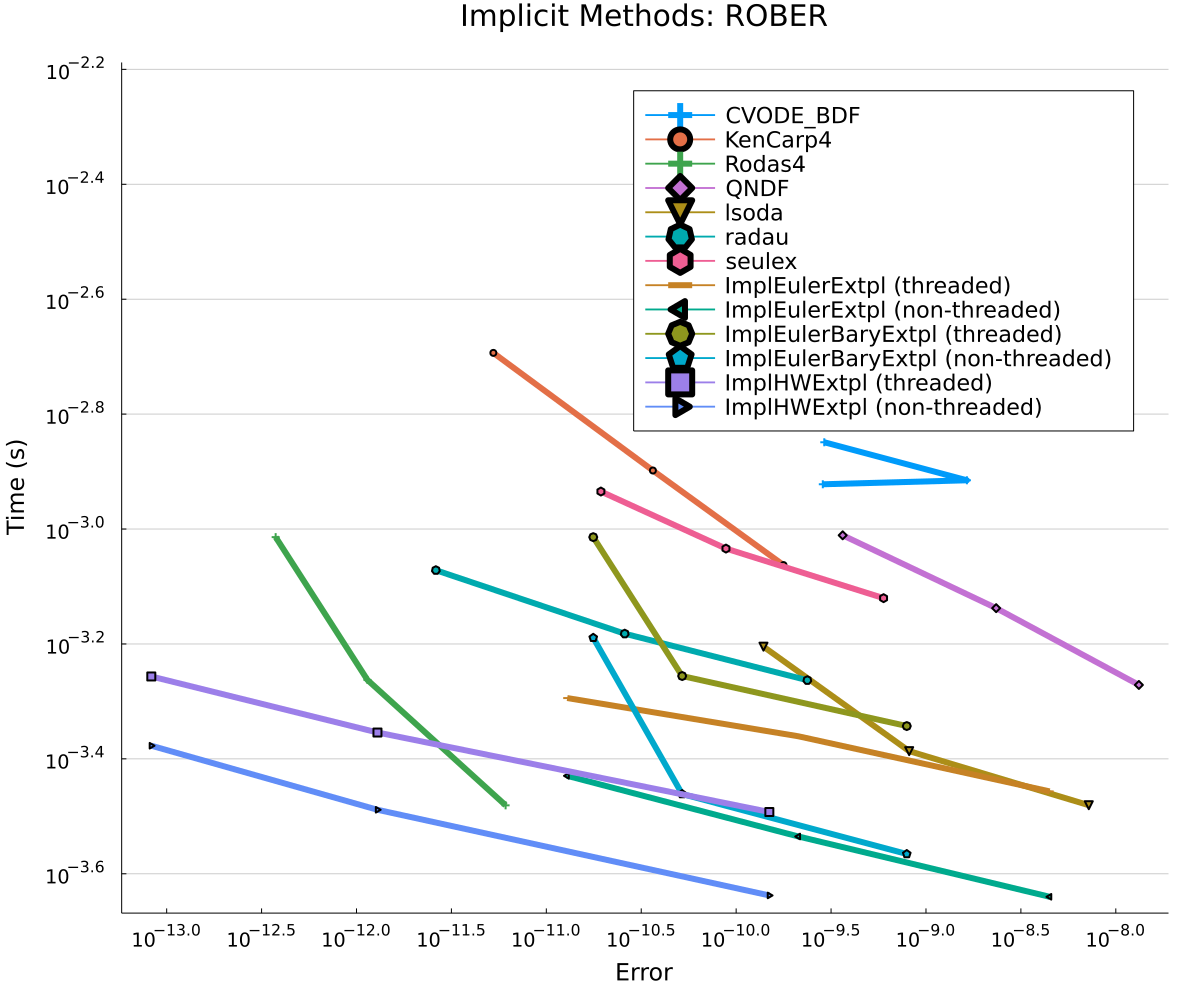}
\label{fig_4}}
\hfil
\vspace{-25px}
\caption{Benchmark on ROBER Problem with low tolerances [\ref{sssec:rober}].}\label{fig:robercomparsion}
\vspace{-15px}
\end{figure}

To establish these methods as efficient versions of extrapolation methods, we benchmarked these new implementations against standard optimized extrapolation implementations. The most widely used optimized implementations of extrapolation methods come from Hairer's FORTRAN suite \cite{hairer2004fortran}, with ODEX \cite{hairer1993solving} as a GBS \cite{bulirsch1966numerical} extrapolation method and SEULEX as a linear implicit Euler method. Figure \ref{fig:linearcomparison} illustrates that our implementation outperforms the ODEX \cite{hairer1993solving} FORTRAN method by 6x on the 100 Independent Linear ODEs problem. And similar to what was shown in Ketcheson 2014 \cite{ketcheson2014comparison}, we see that the parallelized extrapolation methods were able to outperform the dop853 high order Runge-Kutta method implementation of Hairer by around 4x (they saw closer to 2x). However, when we compare the explicit extrapolation methods against newer optimized Runge-Kutta methods like Verner's efficient 9th order method \cite{verner2010numerically}, we only matched the SOA\footnote{In the shown Figure \ref{fig:linearcomparison}, the methods are approximately 1.5x faster, though on this random linear ODE benchmark, different random conditions lead to different performance results, averaging around 0.8x-2x, making the speedup generally insignificant}, which we suspect is due to the explicit methods being able to exploit less parallelism than the implicit variants. Thus for the rest of the benchmarks we focused on the implicit extrapolation methods. In Figure \ref{fig:robercomparsion} we benchmarked the implicit extrapolation methods against the Hairer SEULEX \cite{hairer2004fortran} implementation on the ROBER \cite{robertson1976numerical} problem and demonstrated an average of 4x acceleration over the Hairer implementations, establishing the efficiency of this implementation of the algorithm.

\subsection{State-Of-The-Art Performance for Small Stiff ODE Systems}

To evaluate performance of extrapolation methods, we benchmarked the parallel implicit extrapolation implementations on the following set of standard test problems \cite{hairer1991solving}:

\begin{enumerate}
    \item the Robertson equation \cite{robertson1976numerical}, 3 ODEs
    \item Orego (Oregonator) \cite{zhabotinsky2007belousov, field1974oscillations},  3 ODEs
    \item Hires \cite{schafer1975new}, 8 ODEs
    \item Pollution \cite{verwer1994gauss}, 20 ODEs
    \item QSP \cite{bidkhori2012modeling}, 109 ODEs
\end{enumerate}

Figures 2-6 show the Work-Precision Diagrams with each of the respective problems with low tolerances. Implicit Euler Extrapolation outperforms other solvers with lower times at equivalent errors. However, very small systems (3-7 ODEs) show a disadvantage when threading is enabled, indicating that overhead is not overcome at this size, but the non-threaded version still achieves SOA by about 2x over the next best methods (Rodas4 and radau) and matching lsoda on HIRES. By 20 ODEs, the multi-threaded form becomes more efficient and is the SOA algorithm by roughly 2.5x. But by the 109 ODEs of the QSP model, the multi-threaded form is noticeably more efficient than the non-multithreaded and the extrapolation methods  outperformed lsoda BDF and CVODE's BDF method by about 2x. Together these results show the implicit extrapolation methods as SOA for small ODE systems, with a lower bound on when multi-threading is beneficial.

\begin{figure}[!t]
\centering
{\includegraphics[width=\linewidth]{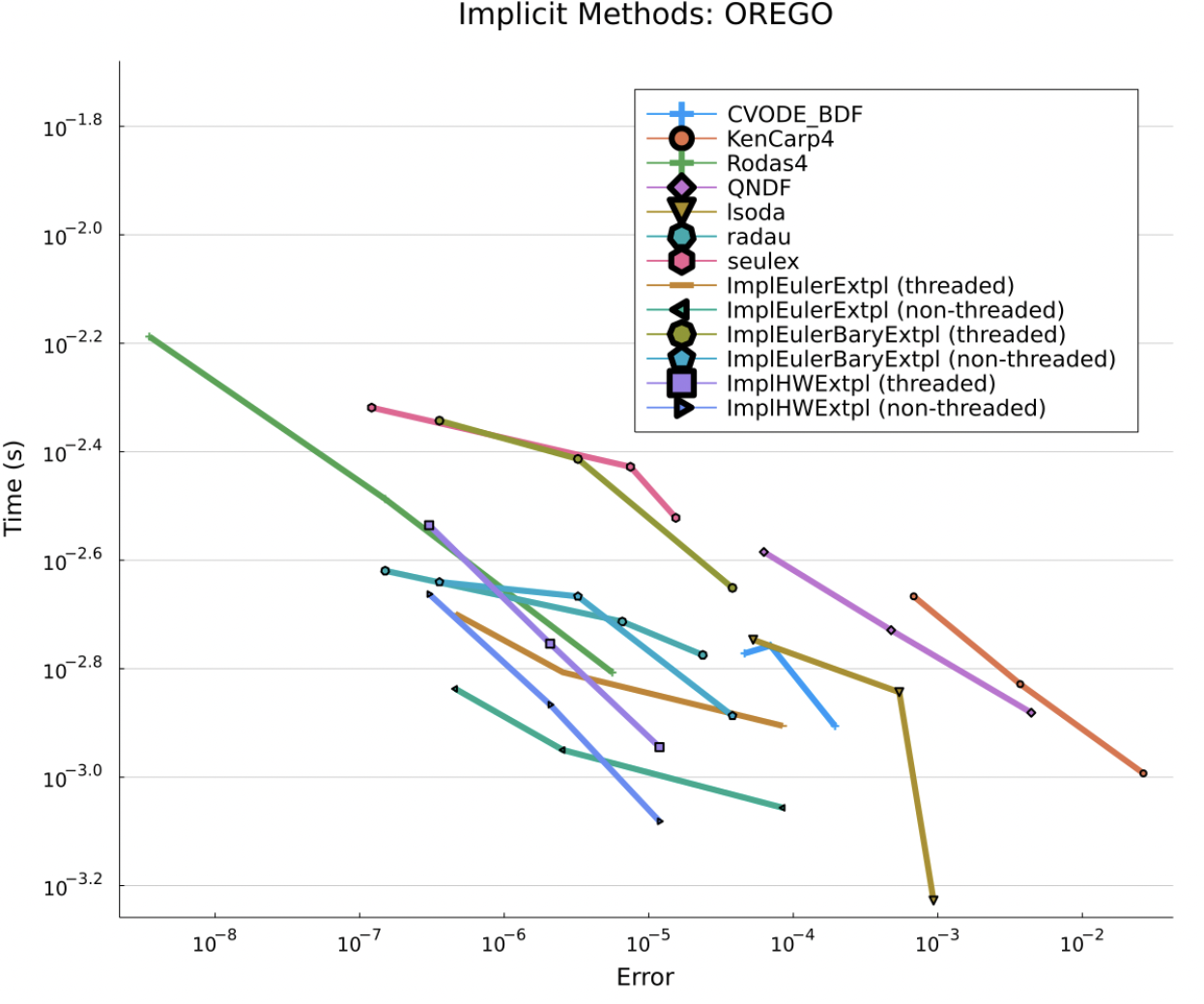}
\label{fig_5}}
\hfil
\vspace{-25px}
\caption{Benchmark on Orego Problem with low tolerances [\ref{sssec:orego}].}
\label{fig:oregocomparsion}
\vspace{-10px}
\end{figure}

\begin{figure}[!t]
\centering
{\includegraphics[width=\linewidth]{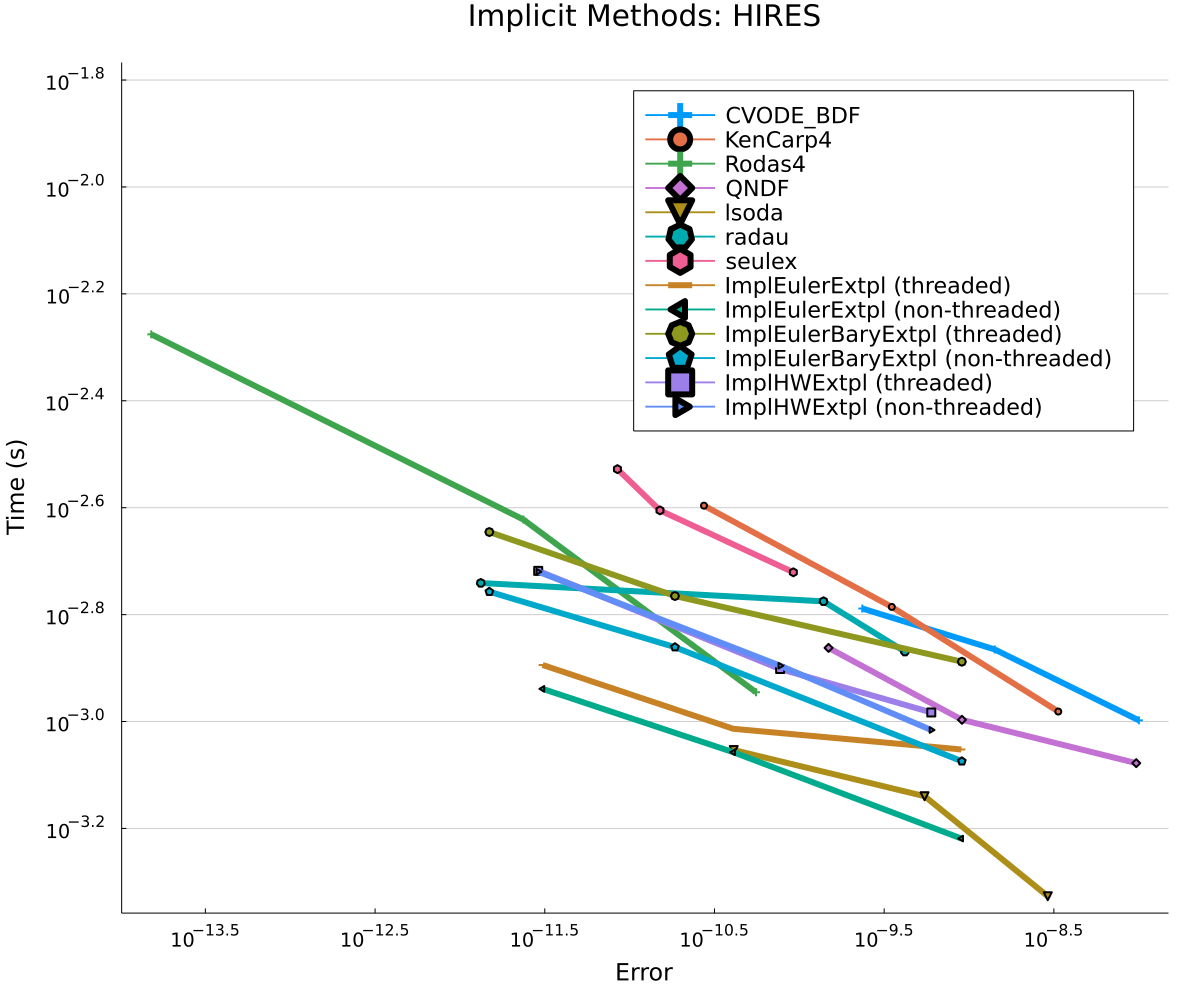}
\label{fig_6}}
\vspace{-25px}
\caption{Benchmark on Hires Problem with low tolerances [\ref{sssec:hires}].}
\label{fig:hirescomparsion}
\vspace{-10px}
\end{figure}

\begin{figure}
\centering
{\includegraphics[width=\linewidth]{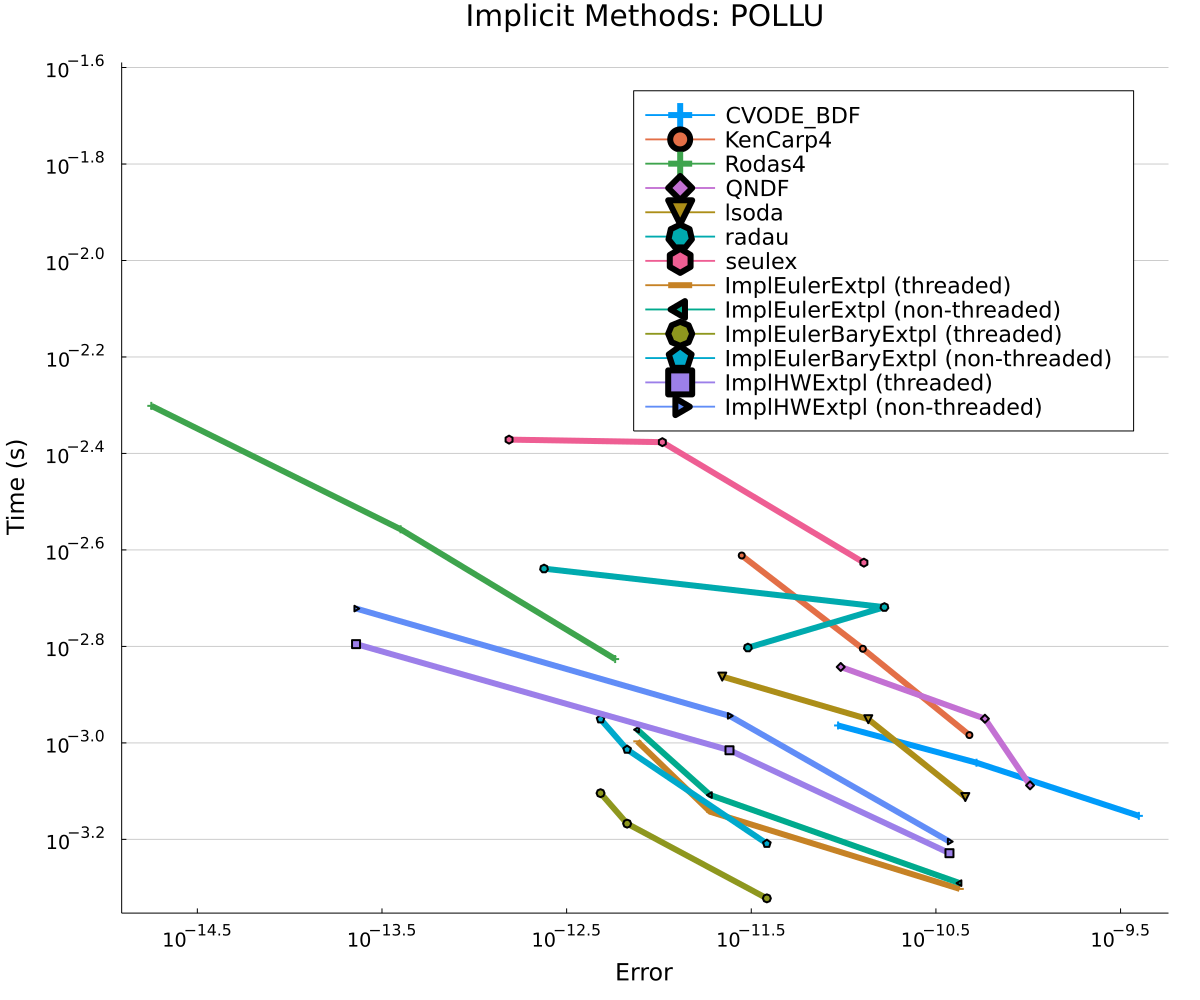}
\label{fig_7}}
\vspace{-25px}
\caption{Benchmark on Pollution Problem with low tolerances [\ref{sssec:pollution}].}
\label{fig:pollucomparsion}
\vspace{-10px}
\end{figure}

\begin{figure}[!t]
\centering
\includegraphics[width=\linewidth]{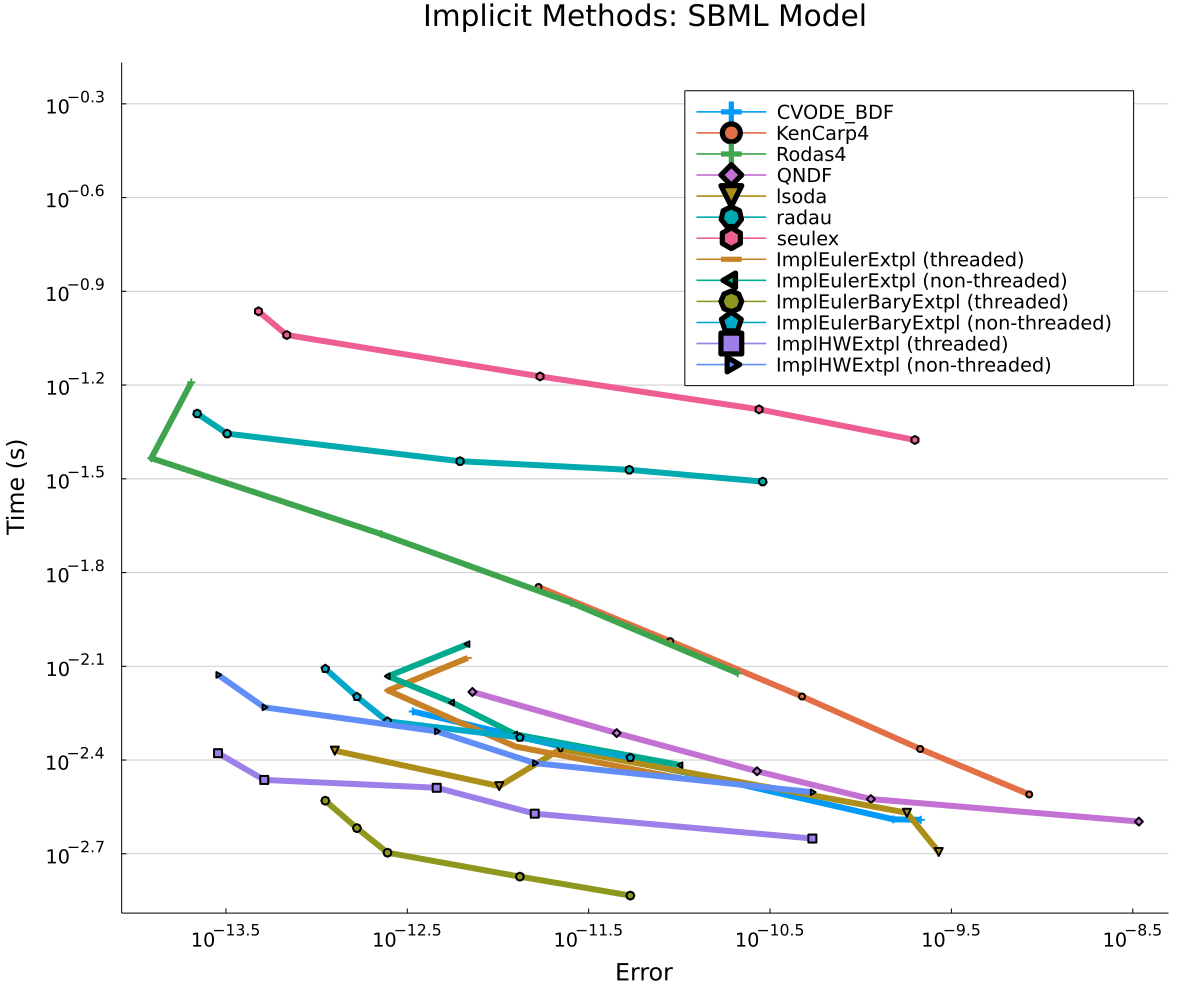}
\vspace{-25px}
\caption{Benchmark on QSP model with low tolerances [\ref{sssec:qspmodel}].}
\label{fig:qspcomparsion}
\vspace{-15px}
\end{figure}

\section{Discussion}

If one follows the tutorial coding examples for the most common differential equation solver suites in high-level languages (R's deSolve \cite{soetaert2010solving}, Python's SciPy \cite{virtanen2020scipy} , Julia's DifferentialEquations.jl \cite{rackauckas2017differentialequations}, MATLAB's ODE Suite \cite{shampine1997matlab}), the vast majority of users in common case of $<200$ ODEs will receive no parallel acceleration even though parallel hardware is readily available\footnote{Unless the $f$ function is sufficiently expensive, like when it is a neural network evaluation}. In this manuscript we have demonstrated how the theory of within-method parallel extrapolation methods can be used to build a method which achieves SOA for the high accuracy approximation regime on this typical ODE case. We distribute this method as open source software: existing users of DifferentialEquations.jl \cite{rackauckas2017differentialequations} can make use of this technique with no changes to their code beyond the characters for the algorithm choice. As such, this is the first demonstration the authors know of for automated acceleration of ODE solves on typical small-scale ODEs using parallel acceleration.\footnote{The codes can be found at:\\ \url{https://github.com/utkarsh530/ParallelExBenchmarks.jl}} 

Lastly, the extrapolation methods were chosen as the foundation because their arbitrary order allows the methods to easily adapt the number of parallel compute chunks to different core counts. This implementation and manuscript statically batched the compute in a way that is independent of the specific core count, leading to generality in the resulting method. However, the downside of this approach is that extrapolation methods are generally inefficient in comparison to other methods such as Rosenbrock \cite{rosenbrock1963some} or SDIRK \cite{alexander1977diagonally} methods. In none of the numerous stiff ODE benchmarks by Harier II \cite{hairer1991solving} or the SciML suite\footnote{https://github.com/SciML/SciMLBenchmarks.jl} do the implicit extrapolation methods achieve SOA performance. Thus while the parallelization done here was able to achieve SOA, even better results could likely be obtained by developing parallel Rosenbrock \cite{rosenbrock1963some} or SDIRK \cite{alexander1977diagonally} methods (using theory described in \cite{burrage1995parallel}) and implementing them using our parallelization techniques. This would have the trade-off of being tied to a specific core count, but could achieve better performance through tableau optimization.

\section*{Acknowledgments}

This material is based upon work supported by the National Science Foundation under grant no. OAC-1835443, SII-2029670,
ECCS-2029670, OAC-2103804, and PHY-2021825.  We also gratefully acknowledge the U.S. Agency for
International Development through Penn State for grant no. S002283-USAID. The information, data, or work presented herein was
funded in part by the Advanced Research Projects Agency-Energy (ARPA-E), U.S. Department of Energy, under Award Number DE-AR0001211
and DE-AR0001222. We also gratefully acknowledge the U.S. Agency for International Development through Penn State for grant no.
S002283-USAID. The views and opinions of authors expressed herein do not necessarily state or reflect those of the United States
Government or any agency thereof. This material was supported by The Research Council of Norway and Equinor ASA through Research
Council project "308817 - Digital wells for optimal production and drainage". Research was sponsored by the United States Air Force
Research Laboratory and the United States Air Force Artificial Intelligence Accelerator and was accomplished under Cooperative
Agreement Number FA8750-19-2-1000. The views and conclusions contained in this document are those of the authors and should not be
interpreted as representing the official policies, either expressed or implied, of the United States Air Force or the U.S. Government.
The U.S. Government is authorized to reproduce and distribute reprints for Government purposes notwithstanding any copyright notation herein.

\bibliographystyle{IEEEtran}
\bibliography{IEEEabrv,refs.bib}

\begin{thebibliography}{10}
\providecommand{\url}[1]{#1}
\csname url@samestyle\endcsname
\providecommand{\newblock}{\relax}
\providecommand{\bibinfo}[2]{#2}
\providecommand{\BIBentrySTDinterwordspacing}{\spaceskip=0pt\relax}
\providecommand{\BIBentryALTinterwordstretchfactor}{4}
\providecommand{\BIBentryALTinterwordspacing}{\spaceskip=\fontdimen2\font plus
\BIBentryALTinterwordstretchfactor\fontdimen3\font minus
  \fontdimen4\font\relax}
\providecommand{\BIBforeignlanguage}[2]{{%
\expandafter\ifx\csname l@#1\endcsname\relax
\typeout{** WARNING: IEEEtran.bst: No hyphenation pattern has been}%
\typeout{** loaded for the language `#1'. Using the pattern for}%
\typeout{** the default language instead.}%
\else
\language=\csname l@#1\endcsname
\fi
#2}}
\providecommand{\BIBdecl}{\relax}
\BIBdecl

\bibitem{dongarra1990set}
J.~J. Dongarra, J.~Du~Croz, S.~Hammarling, and I.~S. Duff, ``A set of level 3
  basic linear algebra subprograms,'' \emph{ACM Transactions on Mathematical
  Software (TOMS)}, vol.~16, no.~1, pp. 1--17, 1990.

\bibitem{anderson1999lapack}
E.~Anderson, Z.~Bai, C.~Bischof, L.~S. Blackford, J.~Demmel, J.~Dongarra,
  J.~Du~Croz, A.~Greenbaum, S.~Hammarling, A.~McKenney \emph{et~al.},
  \emph{LAPACK users' guide}.\hskip 1em plus 0.5em minus 0.4em\relax SIAM,
  1999.

\bibitem{Hindmarsh_2005}
\BIBentryALTinterwordspacing
A.~C. Hindmarsh, P.~N. Brown, K.~E. Grant, S.~L. Lee, R.~Serban, D.~E.
  Shumaker, and C.~S. Woodward, ``{SUNDIALS},'' \emph{{ACM} Transactions on
  Mathematical Software}, vol.~31, no.~3, pp. 363--396, sep 2005. [Online].
  Available: \url{https://doi.org/10.1145%2F1089014.1089020}
\BIBentrySTDinterwordspacing

\bibitem{rackauckas2017differentialequations}
C.~Rackauckas and Q.~Nie, ``Differentialequations.jl--a performant and
  feature-rich ecosystem for solving differential equations in julia,''
  \emph{Journal of Open Research Software}, vol.~5, no.~1, 2017.

\bibitem{bezanson2017julia}
J.~Bezanson, A.~Edelman, S.~Karpinski, and V.~B. Shah, ``Julia: A fresh
  approach to numerical computing,'' \emph{SIAM review}, vol.~59, no.~1, pp.
  65--98, 2017.

\bibitem{burrage1995parallel}
K.~Burrage, \emph{Parallel and sequential methods for ordinary differential
  equations}.\hskip 1em plus 0.5em minus 0.4em\relax Clarendon Press, 1995.

\bibitem{xbraid-package}
``{XB}raid: Parallel multigrid in time,'' \url{http://llnl.gov/casc/xbraid}.

\bibitem{xbraid-tutorial}
J.~Schroder and R.~Falgout, ``Xbraid tutorial,'' in \emph{18th Copper Mountain
  Conference on Multigrid Methods}.\hskip 1em plus 0.5em minus 0.4em\relax
  Colorado, 2017.

\bibitem{nowak1998parallel}
U.~Nowak, R.~Ehrig, and L.~Oeverdieck, ``Parallel extrapolation methods and
  their application in chemical engineering,'' in \emph{International
  Conference on High-Performance Computing and Networking}.\hskip 1em plus
  0.5em minus 0.4em\relax Springer, 1998, pp. 419--428.

\bibitem{davis2004algorithm_umfpack}
T.~A. Davis, ``Algorithm 832: Umfpack v4. 3---an unsymmetric-pattern
  multifrontal method,'' \emph{ACM Transactions on Mathematical Software
  (TOMS)}, vol.~30, no.~2, pp. 196--199, 2004.

\bibitem{davis2014suitesparse}
T.~Davis, W.~Hager, and I.~Duff, ``Suitesparse,'' \emph{URL: faculty. cse.
  tamu. edu/davis/suitesparse. html}, 2014.

\bibitem{petzold1983automatic}
L.~Petzold, ``Automatic selection of methods for solving stiff and nonstiff
  systems of ordinary differential equations,'' \emph{SIAM journal on
  scientific and statistical computing}, vol.~4, no.~1, pp. 136--148, 1983.

\bibitem{ketcheson2014comparison}
D.~Ketcheson and U.~bin Waheed, ``A comparison of high-order explicit
  runge--kutta, extrapolation, and deferred correction methods in serial and
  parallel,'' \emph{Communications in Applied Mathematics and Computational
  Science}, vol.~9, no.~2, pp. 175--200, 2014.

\bibitem{hairer1993solving}
E.~Hairer, S.~P. Norsett, and G.~Wanner, ``Solving ordinary differential
  equations i: Nonsti problems, volume second revised edition,'' 1993.

\bibitem{hairer1991solving}
Hairer and Peters, \emph{Solving Ordinary Differential Equations II}.\hskip 1em
  plus 0.5em minus 0.4em\relax Springer Berlin Heidelberg, 1991.

\bibitem{althaus_2018}
\BIBentryALTinterwordspacing
K.~Althaus, ``Theory and implementation of the adaptive explicit midpoint rule
  including order and stepsize control,'' 2018. [Online]. Available:
  \url{https://github.com/AlthausKonstantin/Extrapolation/blob/master/Bachelor%20Theseis.pdf}
\BIBentrySTDinterwordspacing

\bibitem{aitken1932interpolation}
A.~Aitken, ``On interpolation by iteration of proportional parts, without the
  use of differences,'' \emph{Proceedings of the Edinburgh Mathematical
  Society}, vol.~3, no.~1, pp. 56--76, 1932.

\bibitem{neville1934iterative}
E.~H. Neville, \emph{Iterative interpolation}.\hskip 1em plus 0.5em minus
  0.4em\relax St. Joseph's IS Press, 1934.

\bibitem{deuflhard1983order}
P.~Deuflhard, ``Order and stepsize control in extrapolation methods,''
  \emph{Numerische Mathematik}, vol.~41, no.~3, pp. 399--422, 1983.

\bibitem{romberg1955vereinfachte}
W.~Romberg, ``Vereinfachte numerische integration,'' \emph{Norske Vid. Selsk.
  Forh.}, vol.~28, pp. 30--36, 1955.

\bibitem{bulirsch1966numerical}
R.~Bulirsch and J.~Stoer, ``Numerical treatment of ordinary differential
  equations by extrapolation methods,'' \emph{Numerische Mathematik}, vol.~8,
  no.~1, pp. 1--13, 1966.

\bibitem{stetter1970symmetric}
H.~J. Stetter, ``Symmetric two-step algorithms for ordinary differential
  equations,'' \emph{Computing}, vol.~5, no.~3, pp. 267--280, 1970.

\bibitem{berrut2004barycentric}
J.-P. Berrut and L.~N. Trefethen, ``Barycentric lagrange interpolation,''
  \emph{SIAM review}, vol.~46, no.~3, pp. 501--517, 2004.

\bibitem{webb2012stability}
M.~Webb, L.~N. Trefethen, and P.~Gonnet, ``Stability of barycentric
  interpolation formulas for extrapolation,'' \emph{SIAM Journal on Scientific
  Computing}, vol.~34, no.~6, pp. A3009--A3015, 2012.

\bibitem{dahlquist1963special}
G.~G. Dahlquist, ``A special stability problem for linear multistep methods,''
  \emph{BIT Numerical Mathematics}, vol.~3, no.~1, pp. 27--43, 1963.

\bibitem{bader1983semi}
G.~Bader and P.~Deuflhard, ``A semi-implicit mid-point rule for stiff systems
  of ordinary differential equations,'' \emph{Numerische Mathematik}, vol.~41,
  no.~3, pp. 373--398, 1983.

\bibitem{gragg1964repeated}
W.~B. Gragg, ``Repeated extrapolation to the limit in the numerical solution of
  ordinary differential equations.'' CALIFORNIA UNIV LOS ANGELES, Tech. Rep.,
  1964.

\bibitem{Lindberg_1971}
\BIBentryALTinterwordspacing
B.~Lindberg, ``On smoothing and extrapolation for the trapezoidal rule,''
  \emph{{BIT}}, vol.~11, no.~1, pp. 29--52, mar 1971. [Online]. Available:
  \url{https://doi.org/10.1007%2Fbf01935326}
\BIBentrySTDinterwordspacing

\bibitem{ElrodPolyester}
C.~Elrod, ``Polyester.jl,'' \url{https://github.com/JuliaSIMD/Polyester.jl},
  2021.

\bibitem{hairer2004fortran}
E.~Hairer, ``Fortran and matlab codes,'' 2004.

\bibitem{verner2010numerically}
J.~H. Verner, ``Numerically optimal runge--kutta pairs with interpolants,''
  \emph{Numerical Algorithms}, vol.~53, no.~2, pp. 383--396, 2010.

\bibitem{robertson1976numerical}
H.~Robertson, ``Numerical integration of systems of stiff ordinary differential
  equations with special structure,'' \emph{IMA Journal of Applied
  Mathematics}, vol.~18, no.~2, pp. 249--263, 1976.

\bibitem{zhabotinsky2007belousov}
A.~M. Zhabotinsky, ``Belousov-zhabotinsky reaction,'' \emph{Scholarpedia},
  vol.~2, no.~9, p. 1435, 2007.

\bibitem{field1974oscillations}
R.~J. Field and R.~M. Noyes, ``Oscillations in chemical systems. iv. limit
  cycle behavior in a model of a real chemical reaction,'' \emph{The Journal of
  Chemical Physics}, vol.~60, no.~5, pp. 1877--1884, 1974.

\bibitem{schafer1975new}
E.~Sch{\"a}fer, ``A new approach to explain the “high irradiance responses”
  of photomorphogenesis on the basis of phytochrome,'' \emph{Journal of
  Mathematical Biology}, vol.~2, no.~1, pp. 41--56, 1975.

\bibitem{verwer1994gauss}
J.~G. Verwer, ``Gauss--seidel iteration for stiff odes from chemical
  kinetics,'' \emph{SIAM Journal on Scientific Computing}, vol.~15, no.~5, pp.
  1243--1250, 1994.

\bibitem{bidkhori2012modeling}
G.~Bidkhori, A.~Moeini, and A.~Masoudi-Nejad, ``Modeling of tumor progression
  in nsclc and intrinsic resistance to tki in loss of pten expression,''
  \emph{PloS one}, vol.~7, no.~10, p. e48004, 2012.

\bibitem{soetaert2010solving}
K.~Soetaert, T.~Petzoldt, and R.~W. Setzer, ``Solving differential equations in
  r: package desolve,'' \emph{Journal of statistical software}, vol.~33, pp.
  1--25, 2010.

\bibitem{virtanen2020scipy}
P.~Virtanen, R.~Gommers, T.~E. Oliphant, M.~Haberland, T.~Reddy, D.~Cournapeau,
  E.~Burovski, P.~Peterson, W.~Weckesser, J.~Bright \emph{et~al.}, ``Scipy 1.0:
  fundamental algorithms for scientific computing in python,'' \emph{Nature
  methods}, vol.~17, no.~3, pp. 261--272, 2020.

\bibitem{shampine1997matlab}
L.~F. Shampine and M.~W. Reichelt, ``The matlab ode suite,'' \emph{SIAM journal
  on scientific computing}, vol.~18, no.~1, pp. 1--22, 1997.

\bibitem{rosenbrock1963some}
H.~Rosenbrock, ``Some general implicit processes for the numerical solution of
  differential equations,'' \emph{The Computer Journal}, vol.~5, no.~4, pp.
  329--330, 1963.

\bibitem{alexander1977diagonally}
R.~Alexander, ``Diagonally implicit runge--kutta methods for stiff ode’s,''
  \emph{SIAM Journal on Numerical Analysis}, vol.~14, no.~6, pp. 1006--1021,
  1977.

\bibitem{li2010biomodels}
C.~Li, M.~Donizelli, N.~Rodriguez, H.~Dharuri, L.~Endler, V.~Chelliah, L.~Li,
  E.~He, A.~Henry, M.~I. Stefan \emph{et~al.}, ``Biomodels database: An
  enhanced, curated and annotated resource for published quantitative kinetic
  models,'' \emph{BMC systems biology}, vol.~4, no.~1, pp. 1--14, 2010.

\bibitem{hucka2003systems}
M.~Hucka, A.~Finney, H.~M. Sauro, H.~Bolouri, J.~C. Doyle, H.~Kitano, A.~P.
  Arkin, B.~J. Bornstein, D.~Bray, A.~Cornish-Bowden \emph{et~al.}, ``The
  systems biology markup language (sbml): a medium for representation and
  exchange of biochemical network models,'' \emph{Bioinformatics}, vol.~19,
  no.~4, pp. 524--531, 2003.

\bibitem{ma2021modelingtoolkit}
Y.~Ma, S.~Gowda, R.~Anantharaman, C.~Laughman, V.~Shah, and C.~Rackauckas,
  ``Modelingtoolkit: A composable graph transformation system for
  equation-based modeling,'' 2021.

\end{thebibliography}

\section{Appendix}

\subsection{Models}

The first test problem is the ROBER Problem:
\begin{align}
\begin{split}
{}&\frac{dy_1}{dt} = -k_1y_1 + k_3y_2y_3,
\\
{}&\frac{dy_2}{dt} = k_1y_1 + -k_2y_2^2 -k_3y_2y_3,
\\
{}&\frac{dy_3}{dt} = k_2y_2^2.
\end{split}
\end{align}
The initial conditions are $y = [1.0,0.0,0.0]$ and $k = (0.04,3e7,1e4)$. The time span for integration is $t = (0.0\,s,1e5\,s)$.\cite{robertson1976numerical}, \cite{hairer1991solving}
The second test problem is OREGO:
%
%
\begin{align}
{}&\frac{dy_1}{dt} = -k_1(y_2 + y_1(1-k_2y_1 - y_2)),
\\
{}&\frac{dy_2}{dt} = \frac{y_3 - (1 + y_1)y_2}{k_1},
\\
{}&\frac{dy_3}{dt} = k_3(y_1 - y_3).
\end{align}
The initial conditions are $y = [1.0,2.0,3.0]$ and $k = (77.27,8.375\times10^{-6},0.161)$. The time span for integration is $t = (0.0\,s,30.0\,s)$.
The third test problem is HIRES:
%
%
\begin{align}
\frac{dy_1}{dt} &= -1.71y_1 + 0.43y_2 + 8.32y_3 + 0.0007,
\\
\frac{dy_2}{dt} &= 1.71y_1 - 8.75y_2,
\\
\frac{dy_3}{dt} &= -10.03y_3 + 0.43y_4 + 0.035y_5,
\\ 
\frac{dy_4}{dt} &= 8.32y_2 + 1.71y_3 - 1.12y_4, 
\\ 
\frac{dy_5}{dt} &= -1.745y_5 + 0.43y_6 + 0.43y_7, 
\\ 
\frac{dy_6}{dt} &=  -280.0y_6y_8 + 0.69y_4 \\
                & + 1.71y_5 - 0.43y_6 + 0.69y_7, 
\\
\frac{dy_7}{dt} &= 280.0y_6y_8 - 1.81y_7,
\\
\frac{dy_8}{dt} &= -280.0y_6y_8 + 1.81y_7.
\end{align}
The initial conditions are:
\begin{equation}
 y = [1.0,0.0,0.0,0.0,0.0,0.0,0.0,0.0057].
\end{equation}
The time span for integration is $t = (0.0\,s,321.8122\,s)$.
The fourth test problem is POLLU:
%
%
\begin{align}
\frac{dy_1}{dt} &= -k_1y_1 - k_{10}y_{11}y_1 - k_{14}y_1y_6 - k_{23}y_1y_4 \\
                & - k_{24}y_{19}y_1 + k_2y_2y_4 + k_3y_5y_2 + k_9y_{11}y_2 \\
                & + k_{11}y_{13} + k_{12}y_{10}y_2 + k_{22}y_{19} + k_{25}y_{20},
\\
\frac{dy_2}{dt} &= -k_2y_2y_4 - k_3y_5y_2 - k_9y_{11}y_2 - k_{12}y_{10}y_{2}\\ 
                & + k_1y_1 + k_{21}y_{19},
\\
\frac{dy_3}{dt} &= -k_{15}y_3 + k_1y_1 + k_{17}y_4 + k_{19}y_{16} + k_{22}y_{19},
\\ 
\frac{dy_4}{dt} &= -k_2y_2y_4 - k_{16}y_4 - k_{17}y_4 - k_{23}y_1y_4 + k_{15}y_3, 
\\ 
\frac{dy_5}{dt} &= -k_3y_5y_2 + 2k_4y_7 + k_6y_7y_6 + k_7u_9 \\
                & + k_{13}y_{14} + k_{20}y_{17}y_6, 
\\ 
\frac{dy_6}{dt} &= -k_6y_7y_6 - k_8y_9y_6 - k_{14}y_1y_6 - k_{20}y_{17}y_6 \\
                & + k_3y_5y_2 + 2k_{18}u_{16}, 
\\
\frac{dy_7}{dt} &= -k_4y_7 - k_5y_7 - k_6y_7y_6 + k_{13}y_{14},
\\
\frac{dy_8}{dt} &= k_4y_7 + k_5y_7 + k_6y_7y_6 + k_7y_9,
\\
\frac{dy_9}{dt} &= -k_7y_9 - k_8y_9y_6,
\\
\frac{dy_{10}}{dt} &= -k_{12}y_{10}y_2 + k_7y_9 + k_9y_{11}y_2,
\\
\frac{dy_{11}}{dt} &= -k_9y_{11}y_2 - k_{10}y_{11}y_1 + k_8y_9y_6 + k_{11}y_{13},
\\ 
\frac{dy_{12}}{dt} &= k_9y_{11}y_2, 
\\ 
\frac{dy_{13}}{dt} &= -k_{11}y_{13} + k_{10}y_{11}y_1, 
\\ 
\frac{dy_{14}}{dt} &=  -k_{13}y_{14} + k_{12}y_{10}y_2, 
\\
\frac{dy_{15}}{dt} &= k_{14}y_1y_6,
\\
\frac{dy_{16}}{dt} &= -k_{18}y_{16} - k_{19}y_{16} + k_{16}y_4.
\\
\frac{dy_{17}}{dt} &= -k_{20}y_{17}y_6, 
\\ 
\frac{dy_{18}}{dt} &=  k_{20}y_{17}y_6, 
\\
\frac{dy_{19}}{dt} &= -k_{21}y_{19} - k_{22}y_{19} - k_{24}y_{19}y_1 + k_{23}y_1y_4 + k_{25}y_{20},
\\
\frac{dy_{20}}{dt} &= -k_{25}y_{20} + k_{24}y_{19}y_1.
\end{align}
$\begin{aligned}
k = &[0.35, 26.6, 12300.0, 0.00086, 0.00082, 15000.0, \\
     &0.00013, 24000.0, 16500.0, 9000.0, 0.022, 12000.0, 1.88, \\ 
     &16300.0, 4.8e6, 0.00035, 0.0175, 1.0e8, 4.44e11, \\ 
     &1240.0, 2.1, 5.78, 0.0474, 1780.0, 3.12]. 
\end{aligned}$
$\begin{aligned}
y =&[0.0,0.2,0.0,0.04,0.0,0.0,0.1,0.3,0.017,0.0,\\
    &0.0,0.0,0.0,0.0,0.0,0.0,0.007,0.0,0.0,0.0].
\end{aligned}$

The time span for integration is $t = (0.0\,s,60.0\,s)$.

The QSP model is "BIOMD0000000452 QSP Model, Bidkhori2012 - normal EGFR signalling" \cite{bidkhori2012modeling} from BioModels Database \cite{li2010biomodels}. This is a Systems Biology Markup Model (SBML) \cite{hucka2003systems}, parsed using SBMLToolkit.jl and Catalyst.jl \cite{rackauckas2017differentialequations}, \cite{ma2021modelingtoolkit}. It notably has 109 states and 188 parameters
\subsection{Benchmarks}

\subsubsection{100 Linear ODEs}\label{sssec:linear}

The test tolerances for the solvers in the benchmark Fig. \ref{fig:linearcomparison} are relative tolerances in the range of $(
10^{-7}, 10^{-13})$ and corresponding absolute tolerances are $(10^{-10}, 10^{-16})$.

\begin{table}[htbp]
\caption{Tuned parameters for 100 Linear ODEs}
\begin{center}
\begin{tabular}{|c|c|c|c|}
\hline
Extrapolation Method/Order & minimum& initial & maximum\\
\hline
Midpoint Deuflhard & 5 & 10 & 11\\
\hline
Midpoint Hairer Wanner & 5 & 10 & 11\\
\hline
\end{tabular}
\label{tab:linearparams}
\end{center}
\end{table}

\subsubsection{Rober}\label{sssec:rober}
The test tolerances for the solvers in the benchmark Fig. \ref{fig:robercomparsion} are relative tolerances in the range of $(
10^{-7}, 10^{-9})$
and corresponding absolute tolerances are $(10^{-10}, 10^{-12})$. The integrator used for reference tolerance at $10^{-14}$ is CVODE\_BDF \cite{Hindmarsh_2005}. The parameters for the solvers are tuned to these settings:

\begin{table}[htbp]
\caption{Tuned parameters for ROBER}
\begin{center}
\begin{tabular}{|c|c|c|c|}
\hline
Extrapolation Method/Order & minimum& initial & maximum\\
\hline
Implicit Euler & 3 & 5 & 12\\
\hline
Implicit Euler Barycentric & 4 & 5 & 12\\
\hline
Implicit Hairer Wanner & 2 & 5 & 10\\
\hline
\end{tabular}
\label{tab:roberparams}
\end{center}
\end{table}

\subsubsection{Orego}\label{sssec:orego}
The test tolerances for the solvers in the benchmark Fig. \ref{fig:oregocomparsion} are relative tolerances in the range of $(
10^{-7}, 10^{-9})$
and corresponding absolute tolerances are $(10^{-10}, 10^{-12})$. The integrator used for reference tolerance at $10^{-14}$ is Rodas5 \cite{hairer1991solving}. The parameters for the solvers are tuned to these settings:

\begin{table}[htbp]
\caption{Tuned parameters for OREGO}
\begin{center}
\begin{tabular}{|c|c|c|c|}
\hline
Extrapolation Method/Order & minimum& initial & maximum\\
\hline
Implicit Euler & 3 & 4 & 12\\
\hline
Implicit Euler Barycentric & 3 & 4 & 12\\
\hline
Implicit Hairer Wanner & 2 & 5 & 10\\
\hline
\end{tabular}
\label{tab:oregoparams}
\end{center}
\end{table}

\subsubsection{Hires}\label{sssec:hires}
The test tolerances for the solvers in the benchmark Fig. \ref{fig:hirescomparsion} are relative tolerances in the range of $(
10^{-7}, 10^{-9})$
and corresponding absolute tolerances are $(10^{-10}, 10^{-12})$. The integrator used for reference tolerance at $10^{-14}$ is Rodas5 \cite{hairer1991solving}. The parameters for the solvers are tuned to these settings:

\begin{table}[htbp]
\caption{Tuned parameters for HIRES}
\begin{center}
\begin{tabular}{|c|c|c|c|}
\hline
Extrapolation Method/Order & minimum& initial & maximum\\
\hline
Implicit Euler & 4 & 7 & 12\\
\hline
Implicit Euler Barycentric & 4 & 7 & 12\\
\hline
Implicit Hairer Wanner & 3 & 6 & 10\\
\hline
\end{tabular}
\label{tab:hiresparams}
\end{center}
\end{table}

\subsubsection{POLLUTION}\label{sssec:pollution}
The test tolerances for the solvers in the benchmark Fig. \ref{fig:pollucomparsion} are relative tolerances in the range of $(10^{-8}, 10^{-10})$ and corresponding absolute tolerances are $(10^{-10}, 10^{-13})$. The integrator used for reference tolerance at $10^{-14}$ is CVODE\_BDF \cite{Hindmarsh_2005}. The parameters for the solvers are tuned to these settings:

\begin{table}[htbp]
\caption{Tuned parameters for POLLU}
\begin{center}
\begin{tabular}{|c|c|c|c|}
\hline
Extrapolation Method/Order & minimum& initial & maximum\\
\hline
Implicit Euler & 5 & 6 & 12\\
\hline
Implicit Euler Barycentric & 5 & 6 & 12\\
\hline
Implicit Hairer Wanner & 3 & 6 & 10\\
\hline
\end{tabular}
\label{tab:pollusparams}
\end{center}
\end{table}

\subsubsection{QSP Model}\label{sssec:qspmodel}
The test tolerances for the solvers in the benchmark Fig. \ref{fig:qspcomparsion} are relative tolerances in the range of $(10^{-6}, 10^{-10})$
and corresponding absolute tolerances are $(10^{-9}, 10^{-13})$. The integrator used for reference tolerance at $10^{-14}$ is CVODE\_BDF \cite{Hindmarsh_2005}. The parameters for the solvers are tuned to these settings:

\begin{table}[htbp]
\caption{Tuned parameters for QSP Model}
\begin{center}
\begin{tabular}{|c|c|c|c|}
\hline
Extrapolation Method/Order & minimum& initial & maximum\\
\hline
Implicit Euler & 8 & 9 & 12\\
\hline
Implicit Euler Barycentric & 7 & 8 & 12\\
\hline
Implicit Hairer Wanner & 2 & 5 & 10\\
\hline
\end{tabular}
\label{tab:qspparams}
\end{center}
\end{table}

\end{document}